\documentclass[12pt, leqno]{amsart}

\usepackage{shuffle} 

\usepackage{lscape}
\usepackage{graphicx}
\usepackage{amssymb,amsmath,amsthm,amscd}
\usepackage{mathrsfs}
\usepackage{enumerate}
\usepackage[usenames,dvipsnames]{color}
\usepackage[colorlinks=true, pdfstartview=FitV,
 linkcolor=blue,citecolor=blue,urlcolor=blue]{hyperref}

\usepackage[all]{xy}

\allowdisplaybreaks[4]
\usepackage{oldgerm}
\usepackage{xspace}

\usepackage{marginnote}

\newcommand{\nc}{\newcommand}
\numberwithin{equation}{section}

\newenvironment{red}{\relax\color{red}}{\relax}
\newenvironment{blue}{\relax\color{Dandelion}}{\hspace*{.5ex}\relax}

\newcommand{\beb}{\begin{blue}}
\newcommand{\eb}{\end{blue}}

\newcommand{\berm}[1]{\begin{red}{}\marginnote{\fbox{\scshape\lowercase{M}}}%
#1}  

\newcommand{\berE}[1]{\begin{red}{}\marginnote{\fbox{\scshape\lowercase{E}}}%
#1}  

\newcommand{\berMH}[1]{\begin{red}{}\marginnote{\fbox{\scshape\lowercase{MH}}}%
#1}  

\newcommand{\berS}[1]{\begin{red}{}\marginnote{\fbox{\scshape\lowercase{S}}}%
#1}  

\theoremstyle{plain}
\newtheorem{lemma}{Lemma}[section]
\newtheorem{prop}[lemma]{Proposition}
\newtheorem{theorem}[lemma]{Theorem}

\newcommand{\Prop}{\begin{prop}}
\newcommand{\enprop}{\end{prop}}
\newcommand{\Lemma}{\begin{lemma}}
\newcommand{\enlemma}{\end{lemma}}
\newcommand{\Th}{\begin{theorem}}
\newcommand{\enth}{\end{theorem}}
\newtheorem{corollary}[lemma]{Corollary}
\newcommand{\Cor}{\begin{corollary}}
\newcommand{\encor}{\end{corollary}}
\newtheorem{definition}[lemma]{Definition}
\newtheorem*{conjecture}{Conjecture}
\newcommand{\Def}{\begin{definition}}
\newcommand{\edf}{\end{definition}}
\newtheorem{sublemma}[lemma]{Sublemma}
\newcommand{\Sublemma}{\begin{sublemma}}
\newcommand{\ensub}{\end{sublemma}}

\theoremstyle{definition}
\newtheorem{remark}[lemma]{Remark}

\newtheorem{Convention}[lemma]{Convention}
\newcommand{\Conv}{\begin{Convention}}
\newcommand{\enconv}{\end{Convention}}
\nc{\Conj}{\begin{conjecture}}
\nc{\enconj}{\end{conjecture}}
\nc{\Rem}{\begin{remark}}
\nc{\enrem}{\end{remark}}
\newcommand{\C}{{\mathbb C}}
\newcommand{\Q}{\mathbb {Q}}

\newcommand{\Z}{{\mathbb Z}}
\newcommand{\B}{{\mathbf{B}}}

\newcommand{\D}{\mathscr{D}}

\newcommand{\R}{{\rm R}}

\newcommand{\one}{{\bf{1}}}
\newcommand{\seteq}{\mathbin{:=}}

\newcommand{\hd}{{\operatorname{hd}}}

\newcommand{\g}{{\mathfrak{g}}}

\newcommand{\M}{{\mathscr M}}

\newcommand{\eq}{\begin{eqnarray}}
\newcommand{\eneq}{\end{eqnarray}}

\newcommand{\ms}{\mspace}

\newcommand{\eqn}{\begin{eqnarray*}}
\newcommand{\eneqn}{\end{eqnarray*}}
\newcommand{\on}{\operatorname}

\newcommand{\bni}{\be[{\rm(i)}]}
\newcommand{\bna}{\be[{\rm(a)}]}

\newcommand{\ba}{\begin{array}}
\newcommand{\ea}{\end{array}}

\newcommand{\monoto}{\rightarrowtail}

\newcommand{\eqsub}{\begin{subequations}\begin{eqnarray}}
\newcommand{\eneqsub}{\end{eqnarray}\end{subequations}}

\newcommand{\ol}{\overline}

\nc{\la}{\lambda}
\nc{\lam}{\lambda}
\nc{\U}[1][\g]{U_q(#1)}
\nc{\te}{\tilde{e}}
\nc{\tei}{\tilde{e}_i}
\nc{\tf}{\tilde{f}}
\nc{\tfi}{\tilde{f}_i}
\nc{\tU}{\widetilde U_q(\g)}
\nc{\tE}{\tilde{E}}
\nc{\tF}{\widetilde{\F}}
\nc{\tK}{\widetilde{K}}

\nc{\tk}{\tilde{k}}
\nc{\tkone}{\tk_{\ol{1}}}
\nc{\teone}{\tilde{e}_{\ol{1}}}
\nc{\tfone}{\tilde{f}_{\ol{1}}}

\nc{\teibar}{\tilde{e}_{\ol{i}}} \nc{\tfibar}{\tilde{f}_{\ol{i}}}
\nc{\tki}{{\tk}_{\ol {i}}}

\nc{\BZ}{{\mathbb{Z}}}
\nc{\al}{\alpha}
\nc{\qs}{{q}}
\nc{\lan}{\langle}
\nc{\ran}{\rangle}
\nc{\re}{{\mathrm{re}}}
\nc{\wt}{\operatorname{wt}}
\nc{\ch}{\operatorname{ch}}
\nc{\Um}[1][\g]{U^-_q(#1)}
\nc{\Ue}{U^+_q(\g)}
\nc{\eps}{\varepsilon}
\nc{\vphi}{\varphi}
\nc{\sphi}{\varphi^*}
\nc{\seps}{\varepsilon^*}

\nc{\nn}{\nonumber}

\nc{\vp}{\varpi}
\nc{\cls}{{\operatorname{cl}}}
\nc{\Wt}{{\operatorname{Wt}}}
\nc{\Us}{U'_q(\g)}
\nc{\La}{\Lambda}
\nc{\tLa}{\widetilde\Lambda}
\nc{\ro}{{\rm(}}
\nc{\rf}{{\rm)}}
\nc{\norm}{{\mathrm{norm}}}
\nc{\qbox}{\quad\mbox}
\nc{\braid}{{\mathfrak{B}}}
\nc{\Ad}{\operatorname{Ad}}
\nc{\Aut}{\operatorname{Aut}}
\nc{\dt}[1]{\tilde{\tilde #1}}
\nc{\Sn}{S^{{\mathrm{norm}}}}
\nc{\aff}{{\rm{aff}}}
\nc{\rk}{{\mathrm{rk}}}
\nc{\tP}{\widetilde{P}}
\nc{\tW}{\widetilde{W}}
\nc{\Dyn}{\mathrm{Dyn}}
\nc{\tD}{\widetilde{\Delta}}
\nc{\height}[1]{{\operatorname{ht}}(#1)}
\nc{\bl}{\bigl(}
\nc{\br}{\bigr)}
\nc{\Hecke}{\mathrm{H}}
\nc{\HA}{\Hecke^{\mathrm{A}}}
\nc{\HB}{\Hecke^{\mathrm{B}}}
\newcommand{\scbul}{{\,\raise1pt\hbox{$\scriptscriptstyle\bullet$}\,}}
\nc{\vac}{{\phi}}
\nc{\Bt}{\B_\theta(\g)}
\nc{\be}{\begin{enumerate}}
\nc{\ee}{\end{enumerate}}
\nc{\low}{{\mathrm{low}}}
\nc{\upper}{{\mathrm{up}}}
\nc{\Zodd}{\Z_{\mathrm{odd}}}
\nc{\Ft}[1][n]{\mathbb{P}\mathrm{ol}_{#1}}
\nc{\Ftf}[1][n]{\widetilde{\mathbb{P}\mathrm{ol}}_{#1}}
\nc{\KA}{\on{K}^{\mathrm{A}}}
\nc{\KB}{\on{K}^{\mathrm{B}}}
\nc{\Res}{\on{Res}}
\nc{\Fc}[1][{n,m}]{\mathbf{F}_{#1}}
\nc{\tphi}{\tilde{\varphi}}
\nc{\CO}{\mathscr{O}}
\nc{\inte}{\mathrm{int}}
\nc{\Oint}{\mathcal{O}^{\ge0}_{\inte}}
\nc{\vs}{\vspace*}
\nc{\tLt}{\widetilde{L}}
\nc{\tL}{\widetilde{\Lambda}}
\nc{\tu}{\tilde{u}}
\nc{\noi}{\noindent}
\nc{\heigh}{\mathfrak{t}}
\nc{\lowest}{\mathfrak{l}}
\nc{\rootl}{\mathsf{Q}}
\nc{\cl}{\colon}
\nc{\uqpg}{U'_q(\mathfrak g)}
\nc{\uq}{\uqpg}
\nc{\Oh}{\widehat{\mathcal{O}}}
\nc{\pn}{p_{\mathfrak{n}}}

\nc{\KLR}{KLR algebra}
\nc{\KLRs}{KLR algebras}
\nc{\cor}{\mathbf{k}}
\nc{\cora}{{\cor(A)}}
\nc{\haut}{\mathrm{ht}}
\nc{\tens}{\mathop\otimes}
\nc{\gmod}{\mbox{-$\mathrm{gmod}$}}
\nc{\gMod}{\mbox{-$\mathrm{gMod}$}}
\nc{\proj}{\mbox{-$\mathrm{proj}$}}
\nc{\gproj}{\mbox{-$\mathrm{gproj}$}}
\nc{\smod}{\mbox{-$\mathrm{mod}$}}
\nc{\Mod}{\mbox{-$\mathrm{Mod}$}}
\nc{\h}{\mathfrak h}
\nc{\Rnorm}{R^{\rm{norm}}}
\nc{\Runiv}{R^{\rm{univ}}}
\nc{\Rren}{R^{\rm{ren}}}

\nc{\Vhat}{\widehat{V}}
\nc{\F}{\mathcal{F}}

\def\T{{\mathcal T}}

\nc{\fd}[1][A]{\on{\mathrm{flat.dim}_{#1}}}
\nc{\bP}{{\mathbb{P}}}
\nc{\bPh}{\widehat{\mathbb{P}}}
\nc{\bK}[1][{n}]{\widehat{\mathbb{K}}_{#1}}
\nc{\bV}[1][{n}]{\widehat{V}^{\otimes{#1}}}
\nc{\bVK}[1][{n}]{\widehat{V}^{\otimes{#1}}_{\widehat{\mathbb{K}}}}
\nc{\hV}{\widehat{V}}
\nc{\opp}{\mathrm{opp}}
\nc{\col}{\colon}
\nc{\bnum}{\be[{\rm(i)}]}
\nc{\oep}{\epsilon}
\nc{\qtext}{\quad\text}
\nc{\qtextq}[1]{\quad\text{#1}\quad}
\nc{\longtwoheadrightarrow}[1][]{\xymatrix{\ar@{->>}[r]^-{{#1}}&}}
\nc{\epiTo}[1][]{\longtwoheadrightarrow[{#1}]}
\nc{\epito}{\twoheadrightarrow}
\nc{\monoTo}[1][]{\xymatrix{\ar@{>->}[r]^-{{#1}}&}}
\nc{\sym}{\mathfrak{S}}
\nc{\inp}[1]{{({#1})_{\mathrm{n}}}}
\nc{\rtl}{\rootl}
\nc{\wtd}{\widetilde}
\nc{\etens}{\boxtimes}
\nc{\ds}[1]{\mathrm{d}(#1)}
\nc{\rmat}[1]{{\mathbf{r}}_%
{\mspace{-2mu}\raisebox{-.6ex}{${\scriptstyle{#1}}$}}}
\nc{\rmats}[1]{{\mathbf{r}}_%
{\mspace{-2mu}\raisebox{-.6ex}{${\scriptscriptstyle{#1}}$}}}
\nc{\shc}{\mathcal{C}}
\nc{\shs}{\mathcal{S}}
\nc{\Fct}{{\on{Fct}}}
\nc{\tC}{\widetilde{\shc}}
\nc{\Zp}{\Z_{\ge0}}
\nc{\tPhi}{\widetilde{\Phi}}
\nc{\tT}{{\widetilde{\T}}}
\nc{\Ob}{\on{Ob}}
\nc{\bwr}{\mbox{\large$\wr$}}
\nc{\Img}{\on{Im}}
\nc{\Ab}{\mathcal{A}^{\mathrm{big}}}
\nc{\Sb}{\mathcal{S}^{\mathrm{big}}}
\nc{\As}{\mathcal{A}}
\nc{\Ss}{\mathcal{S}}
\nc{\ntens}{\widetilde{\otimes}}
\nc{\hR}{\widehat{R}}
\nc{\nconv}{\mathop{\mbox{\large $\odot$}}}
\nc{\snconv}{\mbox{\scriptsize$\odot$}}
\nc{\ts}{\tilde{s}}
\nc{\sho}{\mathcal{O}}
\nc{\bc}{\begin{cases}}
\nc{\ec}{\end{cases}}
\nc{\slnh}{{\widehat{\mathfrak{sl}}_N}}
\nc{\UA}{U_q'(\slnh)}
\nc{\KR}{R_K}
\nc{\cQ}{\mathcal{Q}}
\nc{\Irr}{\mathcal{I}rr}
\nc{\tQ}{\widetilde{\cQ}}
\nc{\bs}{\mathbf{s}}
\nc{\bL}{\mathbb{L}}
\nc{\tg}{\tilde{g}}

\nc{\conv}{\mathbin{\mbox{\large $\circ$}}}
\nc{\shconv}{\mathbin{\large\diamond}}
\nc{\sconv}{\mathbin{\large\Delta}}
\nc{\hconv}{\mathbin{\nabla}}

\nc{\Rm}{R^{\mathrm{ren}}}

\nc{\bQ}{\ol{Q}}

\nc{\de}{\on{\textfrak{d}\ms{1mu}}}

\nc{\xmono}{\ar@{>->}}
\nc{\xepi}{\ar@{->>}}
\nc{\db}[1]{\raisebox{-.5ex}[2ex][1.8ex]{$#1$}}
\nc{\wb}[1]{\mbox{$\rule[-1.1ex]{0ex}{2ex}#1$}}
\nc{\univ}{\mathrm{univ}}
\nc{\rM}{{}^*\mspace{-2mu}M}
\nc{\lM}{M^*}
\nc{\uqm}{\uq\smod}
\nc{\tR}{\widetilde{R}_{\gamma,\beta}}
\nc{\tx}{\tilde{x}}
\nc{\bi}{\mathbf{i}}
\nc{\ttau}{\widetilde{\tau}}

\nc{\tEnd}{\on{\widetilde{E}nd}}
\nc{\tHom}{\on{\widetilde{H}om}}

\nc{\K}{{J}}
\nc{\Kex}{{\K}_{\mathrm{ex}}}
\nc{\Kfr}{{\K}_{\mathrm{f\mspace{.01mu}r}}}
\nc{\coro}{\cor}
\nc{\tB}{\widetilde{B}}
\nc{\seed}{\mathscr{S}}

\nc{\up}{\mathrm{up}}
\nc{\bfa}{\mathbf{a}}

\newlength{\mylength}
\nc{\ov}[1]{\overline{#1}}
\nc{\Wlmj}[3]{\W_{#2,#3}^{(#1)}}
\nc{\Mkl}[2]{\M_\ttww(#1,#2)}
\nc{\mqs}{(-q^2)}
\nc{\Cquiver}{\upsigma}

\nc{\mut}[1]{{\mu}_{\mspace{-2mu}\raisebox{-.5ex}{${\scriptstyle{#1}}$}}}

\nc{\Kt}{\mathcal K_t}
\nc{\KT}{\mathbb{K}_t}
\nc{\yim}{y_{i,m}}
\nc{\yjm}{y_{j,m}}
\nc{\yjp}{y_{j,p}}
\nc{\yimp}{y_{i,m+1}}
\nc{\yjmp}{y_{j,m+1}}

\nc{\Refl}{\mathscr{S}}
\nc{\Reflinv}{{\Refl}^{-1}}

\nc{\catC}{\mathscr C}
\nc{\catA}{\mathcal A}
\nc{\shift}{{\mathrm T}}
\nc{\rE}{ \mathsf{E} }
\nc{\rW}{ \mathcal{W} }
\nc{\rES}{ \mathcal{E} }

\nc{\brd}{\sigma} 
\nc{\into}{\xymatrix@C=3ex{{}\ar@{^{(}->}[r]&{}}}
\nc{\dual}{\D}
\nc{\cat}[1][{\g}]{\catC_{#1}^0}
\nc{\qt}[1]{[{#1}]_t}

\nc{\catCO}{{\catC_\g^0}}
\nc{\catCQ}{{\catC_{\qQ}}}
\nc{\Li}{{\La^\infty}}
\nc{\sigZ}{{\sigma_0(\g)}}
\nc{\sigQ}{{\sigma_\qQ(\g)}}

\nc{\ZZ}{{\mathbf{Z}}}
\nc{\sP}{{\mathsf{P}}}
\nc{\sV}{{\mathsf{V}}}
\nc{\rxw}{{\underline{w_0}}}
\nc{\boten}[1]{\overrightarrow{\bigotimes_{#1}}}
\nc{\cmA}{{\mathsf{A}}}
\nc{\cmC}{{\mathsf{C}}}
\nc{\ddD}{{\mathcal{D}}}
\nc{\qQ}{{\mathcal{Q}}}
\nc{\gf}{{\g_{\rm fin}}}
\nc{\If}{{I_{\rm fin}}}
\nc{\cmAf}{{\cmA_{\rm fin}}}
\nc{\weyl}{{\mathsf{W}}}
\nc{\weylf}{{\weyl_{\rm fin}}}
\nc{\Deg}{\mathrm{Deg}}
\nc{\KRc}{{K_{q=1}(R_\cmC\gmod)}}
\nc{\prD}{{\Delta^+_{0}}}
\nc{\n}{{\mathfrak{n}}}
\nc{\Rt}{L} 
\nc{\Cp}{V} 
\nc{\cuspS}{{\mathsf{S}}}
\nc{\st}[1]{\{{#1}\}}
\nc{\WS}{quantum affine Weyl-Schur duality\xspace}
\nc{\CWS}{Quantum affine Weyl-Schur duality}

\title
[{\scalebox{.9}{PBW theoretic approach to the module category of quantum affine algebras}}]{PBW theoretic approach to the module category of quantum affine algebras}

\author[M. Kashiwara]{Masaki Kashiwara}
\thanks{The research of M.\ Kashiwara
was supported by Grant-in-Aid for Scientific Research (B)
15H03608, Japan Society for the Promotion of Science.}
\address[M. Kashiwara]{
Kyoto University Institute for Advanced Study,
Research Institute for Mathematical Sciences, Kyoto University,
Kyoto 606-8502, Japan \& Korea Institute for Advanced Study, Seoul 02455, Korea }
\email[M. Kashiwara]{masaki@kurims.kyoto-u.ac.jp}

\author[M. Kim]{Myungho Kim}
\address[M. Kim]{Department of Mathematics, Kyung Hee University, Seoul 02447, Korea}
\email[M. Kim]{mkim@khu.ac.kr}
\thanks{The research of M.\ Kim was supported by the National Research Foundation of
Korea(NRF) Grant funded by the Korea government(MSIP) (NRF-2017R1C1B2007824).}

\author[S.-j. Oh]{Se-jin Oh}
\thanks{ The research of S.-j.\ Oh was supported by the Ministry of Education of the Republic of Korea and the National Research Foundation of Korea (NRF-2019R1A2C4069647).}
\address[S.-j. Oh]{Department of Mathematics, Ewha Womans University, Seoul 03760, Korea}
\email[S.-j. Oh]{sejin092@gmail.com}

\author[E. Park]{Euiyong Park}

\address[E. Park]{Department of Mathematics, University of Seoul, Seoul 02504, Korea}
\email[E. Park]{epark@uos.ac.kr}


\date{May 10. 2020}

\begin{document}
\maketitle
\begin{abstract}
Let $U_q'(\g)$ be a quantum affine algebra of  untwisted affine 
ADE type and let $\catCO$ be Hernandez-Leclerc's category. 
For a duality datum $\ddD$ in $\catCO$, we denote by $\F_\ddD$ the 
\WS functor.
We give sufficient conditions for a duality datum $\ddD$ to provide the functor $\F_\ddD$ sending simple modules to simple modules. 
Moreover, under the same condition,
the functor $\F_\ddD$ has compatibility 
with the new invariants introduced by the authors.
Then we introduce the notion of cuspidal modules in $\catCO$, and 
show that all simple modules in $\catCO$ can be constructed as the heads of ordered tensor products of cuspidal modules.
We next state that the ordered tensor products of cuspidal modules 
have the unitriangularity property.
\end{abstract}


\section{Introduction}

Let $q$ be an indeterminate and let $\catC_\g$ be the category of finite-dimensional integrable modules over  a quantum affine algebra $U_q'(\g)$.
The category $\catC_\g$ occupies an important position in the  representation
theory of quantum affine algebras because of its rich structure. 
The simple modules in $\catC_\g$ are parameterized using $n$-tuples of polynomials with constant term 1 (called \emph{Drinfeld polynomials}),
which was proved in \cite{CP91,CP94,CP95} for the untwisted cases and in \cite{CP98} for the twisted cases. 
Any simple module can be obtained as the head of an 
ordered tensor product of \emph{fundamental modules}, which was shown in \cite{AK97, Kas02, VV02}.  
A geometric approach to simple modules was also studied in \cite{Nak01, VV02}.

Let $\g_0$ be a finite-dimensional simple Lie algebra of ADE type and $\g$ the untwisted affine Lie algebra associated with $\g_0$.  Hernandez and Leclerc defined the full subcategory $\catCO$ of $\catC_\g$ such that all simple subquotients of its objects appear in tensor products of certain fundamental representations (\cite{HL10}).
Because any simple module in $\catC_\g$ can be obtained as a
tensor product of  suitable parameter shifts of simple modules in $\catCO$, 
the category $\catCO$ contains an essential information of $\catC_\g$.
For each \emph{Q-data}  $\qQ = (Q,\phi)$ of $\g_0$, which is a pair of a Dynkin quiver $Q$ of $\g_0$ and its height function $\phi$, 
Hernandez and Leclerc defined a monoidal subcategory $\catCQ$ of $\catC_\g^0$ and
proved that its complexified Grothendieck ring $\C \otimes _\Z K(\catCQ)$ is isomorphic to the coordinate ring $\C[N]$ of the unipotent group associated with $\g_0$ (\cite{HL15}).
Under this isomorphism, the set of isomorphism classes of simple modules in $\catCQ$ corresponds to the upper global base of $\C[N]$.

The \emph{\WS}  was introduced in \cite{KKK18A}. 
This duality tells us that, for a \emph{duality datum} $\ddD = \{ \Rt_i \}_{i\in J} \subset \catC_\g$, there exists a monoidal functor $\F_\ddD$ from the finite-dimensional graded module category $R\gmod$  of the \emph{quiver Hecke algebra} $R$ (\cite{KL09, R08}) determined by $\ddD$ to the category $\catC_\g$.
In general, it is hard to find conditions for $\ddD$ to provide the functor $\F_\ddD$ with good properties. 
But, to each choice of Q-data $\qQ$, we can associate a \WS $\F_\qQ$ with good properties (\cite{KKK15B, Fu18}): 
$$
\F_\qQ \col R_{\g_0}\gmod \longrightarrow \catCQ \subset \catC_\g^0,
$$
where $R_{\g_0}$ is the symmetric quiver Hecke algebra associated with $\g_0$. 
This functor $\F_\qQ$ sends simple modules of $R_{\g_0}\gmod$ to simple modules of $\catCQ$ and gives an isomorphism at the Grothendieck ring level.

\medskip
One of the main results of this paper is to describe 
 a sufficient condition 
for a duality datum $\ddD$ to provide the functor $\F_\ddD$ with good properties. 
Let $U_q'(\g)$ be an \emph{arbitrary} quantum affine algebra. 
Let $\ddD = \{\Rt_i \}_{i\in J} \subset \catC_\g $ be a duality datum associated with a simply-laced finite Cartan matrix $\cmC$ and let $R_\cmC$ be the symmetric quiver Hecke algebra associated with $\cmC$. 
If $L_i$ are root modules (see \eqref{eq:root})
and $\ddD$ satisfies condition $\eqref{Eq: main thm condition}$ below, 
we say that $\ddD$ is a {\em strong duality datum}.
We show that, for a strong duality datum $\ddD$,
the duality functor $\F_\ddD$ sends simple modules to simple modules and it induces an injective ring homomorphism from $\KRc$ to $K(\catC_\g  )$ 
(see Theorem \ref{Thm: Main1}).
Moreover,
the duality functor is compatible with the new invariants:
\eq\label{eq:pr}
&\phantom{aaaaa}&\ba{l}
\La(M,N) = \La( \F_\ddD(M), \F_\ddD(N) ),\\
\de(M,N) = \de( \F_\ddD(M), \F_\ddD(N) ),\\
( \wt(M), \wt(N)) = - \Li( \F_\ddD(M), \F_\ddD(N) )
\ea
\eneq
for any simple modules $M$, $N$ in $ R_\cmC\gmod$.
Here, $\La $, $\de$ and $\Li $
are new invariants  for 
pairs of objects of $\catC_\g $ introduced in \cite{KKOP19C}.
These invariants are quantum affine algebra analogues of the invariants (with the same notations) for the quiver Hecke algebras.
Note that the block decompositions for $\catC_\g$ and $\catCO$ were given by using the new invariant $\Li$ in \cite{KKOP20}.

\medskip
The other main theorem of this paper is to construct all simple modules in $\catCO$ as the heads of ordered tensor products of \emph{cuspidal modules},
which can be understood as a generalization of the simple module construction using ordered tensor products of fundamental modules
(\cite{CP95, Nak04}).

Suppose that $\g_0$ is a finite-dimensional simple Lie algebra of $ADE$ type 
and $\g$ is the untwisted affine Lie algebra associated with $\g_0$. 
Let $\qQ$ be a Q-data of $\g_0$ and let $\ddD_\qQ$ be the duality datum induced by $\qQ$. Let $\F_\qQ\seteq\F_{\ddD_\qQ}$ be the duality functor associated with
$\ddD_\qQ$. 
Since  $\ddD_\qQ$ is a strong duality datum, 
$\F_\qQ$ sends simples to simples and preserves the new invariants, i.e., \eqref{eq:pr} holds.

Let $w_0$ be  the longest element
of the Weyl group $\weyl_0$ of $\g_0$, and $\ell$ the length of $w_0$.
We choose an arbitrary reduced expression $\rxw$ of $w_0$. 
We define the cuspidal modules  $\{ \cuspS_k \}_{ k\in \Z } \subset \catCO $ to be the simple $U_q'(\g)$-modules given by 
\bna
\item $\cuspS_k = \F_\qQ(\Cp_k)$ for any $k=1, \ldots, \ell$,
\item $\cuspS_{k+\ell} = \dual( \cuspS_k )$ for any $k\in \Z$,
\ee
where $\{\Cp_k \}_{k=1, \ldots, \ell} \subset R_{\g_0}\gmod$ are  the cuspidal modules associated with $\rxw$, and 
$\dual(X)$ and $\dual^{-1}(X)$ denote the right dual and the left dual
of a module $X\in \catC_\g$, respectively.
Note that the cuspidal module $\Cp_k$ corresponds to the \emph{dual PBW vectors} associated with $\rxw$ under the categorification using quiver Hecke algebras.

 For any $\bfa = (a_k)_{k\in \Z} \in \ZZ\seteq\Z_{\ge0}^{\oplus \Z}$, we define 
the ordered tensor product by
$$
\sP_{Q, \rxw} (\bfa) \seteq 
\cdots\tens\cuspS_2^{\otimes a_2} \tens\cuspS_1^{\otimes a_1} \otimes  \cuspS_{0}^{\otimes a_{0}}
 \otimes \cuspS_{-1}^{\otimes a_{-1}}\tens\cuspS_{-2}^{\otimes a_{-2}}\tens \cdots. 
$$
We prove that the head $ \hd \bl\sP_{Q, \rxw} (\bfa)\br$ of the ordered tensor product $\sP_{Q, \rxw} (\bfa)$ is simple and
the simple module $ \hd \bl\sP_{Q, \rxw} (\bfa)\br$
appears only once in $\sP_{Q, \rxw} (\bfa)$.
Moreover, for any simple module $V \in \catCO$, there exists a unique $\bfa\in \ZZ$ such that $V$ is isomorphic to the head $\hd \bl\sP_{Q, \rxw} (\bfa)\br$.
Such an $\bfa$ is denoted by $\bfa_{Q, \rxw}(V)$.

Thus, setting $ \sV _{Q, \rxw} (\bfa) \seteq \hd \bl\sP_{Q, \rxw} (\bfa)\br$,
the set $\{   \sV _{Q, \rxw} (\bfa) \mid \bfa \in \ZZ \}$ is a complete and irredundant set of simple modules of $\catCO$ up to isomorphisms (see Theorem \ref{Thm: main pbw}).
We prove further that,
if $V$ is a simple subquotient of $\sP_{Q, \rxw} (\bfa)$ which is not isomorphic to $ \sV _{Q, \rxw} (\bfa) $, then  
$$
\bfa_{Q, \rxw}(V) \prec \bfa,
$$
where $\prec$ is the bi-lexicographic order on $\ZZ$ (see \eqref{eq:bilexico}).
Thus the family of modules $\st{\sP_{Q, \rxw} (\bfa)}_{\bfa\in\ZZ}$ has the unitriangularity property with respect to $\prec$ (see Theorem \ref{Thm: main tri}).

We can generalize the above results to an arbitrary
quantum affine algebra $U_q'(\g)$ including \emph{twisted type} by using certain strong duality datum in $\catC_\g$ (see Remark~\ref{rmk: generlization of section PBW}).

\smallskip
This paper is an announcement whose details will appear elsewhere.

\vskip 1em 

\section{\WS} \

Let $\cor$ be the algebraic closure of the subfield $\Q(q)$ in the algebraically closed field $\bigcup_{m>0} \C( (q^{1/m}) )$.
Let $U_q'(\g)$ be the  quantum affine algebra over the base field $\cor$ associated with an affine Cartan matrix $\cmA = (a_{i,j})_{i,j\in I}$, and set
 $\catC_\g$ to be the category of finite-dimensional integrable $U_q'(\g)$-modules. 

Let $\cmC = (c_{i,j})_{i,j\in J}$ be  a simply-laced finite Cartan matrix and $R_\cmC$ the  symmetric quiver Hecke algebra associated with $\cmC$.
We denote by $K(R_\cmC\gmod)$ the Grothendieck ring of the category $\R_\cmC\gmod$ of finite-dimensional graded $R_\cmC$-modules.
Note that $K(R_\cmC\gmod)$ is isomorphic to the quantum unipotent coordinate ring $A_{q}(\cmC)_{\Z[q,q^{-1}]}$ (\cite{KL09, R08}) 
and the set of isomorphism classes of simple $R_\cmC$-modules corresponds to the upper global basis of $A_{q}(\cmC)_{\Z[q,q^{-1}]}$ (\cite{R11, VV09}).  
We set $\KRc$ to be the specialization of $K(R_\cmC\gmod)$ at $q=1$.

We freely use new invariants $\La, \Li$ and $\de$ for pairs
of modules in $\catC_\g$  introduced in \cite{KKOP19C}.

A simple $U_q'(\g)$-module $M$ is \emph{real} if $M\otimes M$ is simple.
A {\em root module} is a real simple module $L$ such that
\eq\label{eq:root}
&\phantom{aaaa}&\de(L, \dual^k(L)) = 
\delta(k=\pm 1)
\qtext{for any $k\in\Z$.}
\eneq

Let $\ddD \seteq \{ \Rt_i \}_{i\in J} \subset \catC_\g$
be a family of simple modules of $\catC_\g$.
The family $\ddD $ is called a \emph{duality datum} associated with $\cmC$ if it satisfies the following:
\bna
\item for each $i\in J$, $\Rt_i$ is a real simple module,
\item for any $i,j\in J$ with $i\ne j$, $\de(\Rt_i, \Rt_j) = -c_{i,j}$.
\ee
Then one can construct a monoidal functor 
$$
\F_\ddD \col R_\cmC\gmod \longrightarrow \catC_\g
$$
 using the duality datum $\ddD$ (see \cite{KKK18A, KP}).
Moreover, $\F_\ddD$ is an exact functor.
The functor $\F_\ddD$ is called a \emph{\WS functor} or shortly a \emph{duality functor}.

\smallskip 

A {\em strong duality datum}
is a duality datum $ \ddD= \{ \Rt_i \}_{i\in J}$ such that
all $\Rt_i$ are root modules and
\eq \label{Eq: main thm condition}
&&\phantom{aaaaa}\de(\Rt_i, \dual^k(\Rt_j)) = 
- \delta(k= 0) c_{i,j}
\eneq
for any $k\in \Z$ and $i,j\in J$ such that $i\not=j$.

\begin{theorem} \label{Thm: Main1}
Let $ \ddD= \{ \Rt_i \}_{i\in J}$ be a {\em strong} duality datum associated with a simply-laced finite Cartan matrix $\cmC = (c_{i,j})_{i,j\in J}$.
Then we obtain the following.
\bni
\item The duality functor $\F_\ddD$ sends simple modules to simple modules.
\item The duality functor $\F_\ddD$ induces an injective ring homomorphism 
$$\KRc\monoto K(\catC_\g  ).$$
\ee

\end{theorem}

\smallskip 

The duality functor also has compatibility with the new invariants.

\smallskip 

\begin{theorem} \label{Thm: main 2}
Let $ \ddD= \{ \Rt_i \}_{i\in J}$ be a strong duality datum. 
Then, for any simple modules $M$, $N$ in $ R_\cmC\gmod$, we have 
\bnum
\item $\La(M,N) = \La( \F_\ddD(M), \F_\ddD(N) )$,
\item $\de(M,N) = \de( \F_\ddD(M), \F_\ddD(N) )$,
\item $(\wt M,\wt N) =-\Li( \F_\ddD(M), \F_\ddD(N) )$,
\item $\de\bl \dual^k\F_\ddD(M), \F_\ddD(N)\br=0$ for any $k\not=0,\pm1$,
\item \raisebox{-.7ex}{$\ba[t]{rl}\tL(M,N)&=\de\bl \dual\F_\ddD(M),\F_\ddD(N)\br\\
&=\de\bl \F_\ddD(M), \dual^{-1}\F_\ddD(N)\br.\ea$}
\ee
\end{theorem}

\vskip 1em 

\section{ PBW theoretic approach} \label{sec: PBW}

Let $\g_0$ be a finite-dimensional simple Lie algebra of ADE type,
and $\g$ the \emph{untwisted} affine Kac-Moody algebra associated with $\g_0$. 
We denote by  $\cmA=(a_{i,j})_{i,j\in I}$ the affine Cartan matrix of $\g$ and by $I_0 \subset I$ the index set corresponding to $\g_0$.
For each $i\in I_0$, let $V(\varpi_i)$ be the $i$-th \emph{fundamental representation}.  

Following \cite{HL10},  we denote by $\catCO$ the smallest full subcategory of the category $\catC_\g$ such that 
\bna
\item it contains $\{ V(\varpi_i)_{(-q)^p} \mid i\in I_0,\ p \equiv  d(1,i) \mod2 \}$, 
\item it is stable under taking subquotients, extensions, and tensor products.
\ee
Here $ 1 $ is an arbitrarily chosen element of $I_0$ and $d(1, i)$ is the distance between $1$ and $i$ in the Dynkin diagram of $\g_0$.

A \emph{Q-data} $\qQ = ( Q, \phi_Q  )$ is a pair of a Dynkin quiver $Q$ of $\g_0$
 and a function $\phi_Q(i) \col I_0 \rightarrow \Z$ such that $\phi_Q(1) \in 2\Z$ and $\phi_Q(i) = \phi_Q(j)+1$ for any arrow $i\rightarrow j$ in $Q$.
The function $\phi_Q$ is called a \emph{height function} of $Q$.
For any Q-data $\qQ$, we denote by $\catCQ$ the full monoidal subcategory of $\catCO$ introduced in \cite{HL15}.
For $i\in I_0$, let $\Rt_i$ be the fundamental representation in $\catCQ$ corresponding to the simple root $\al_i$ in the Auslander-Reiten quiver of $\qQ$. 
Then $\{ \Rt_i\}_{i\in I_0}$ forms a strong duality datum
and it induces a duality functor 
$$
\F_\qQ\col R_{\g_0}\gmod \longrightarrow \catCQ \subset \catCO
$$ 
(\cite{KKK15B, Fu18}).
Here $R_{\g_0}$ is the symmetric quiver Hecke algebra associated with $\g_0$. 
Note that the functor $ \F_\qQ $ is an equivalence of categories between $R_{\g_0}\gmod$ and $\catCQ$ after forgetting grading (\cite{Fu18}).

\smallskip 
Hence Theorem~\ref{Thm: main 2} implies the following result.
\begin{theorem} \label{Thm: main3}
Let $\qQ$ be a Q-data of $\g_0$. For any simple modules $M$ and $N$ in $ R_{\g_0}\gmod$,
 we have 
\bni
\item $\La(M,N) = \La( \F_\qQ(M), \F_\qQ(N) )$, 
\item $ ( \wt(M), \wt(N)) = - \Li( \F_\qQ(M), \F_\qQ(N) )$.
\ee
\end{theorem}

\smallskip

Let $\prD$ be the set of positive roots of $\g_0$ and let $\weyl_0 = \langle r_i \mid i\in I_0 \rangle$ be the Weyl group associated with $\g_0$, where $r_i$ is the $i$-th reflection. 
Let $w_0$ be the longest element of $\weyl_0$, and $\ell$ denotes the length of $w_0$. 
We choose an arbitrary reduced expression $\rxw = r_{i_1} r_{i_2} \cdots r_{i_\ell}$ of the longest element $w_0$ of $\weyl_0$.
(We do not assume that $\rxw$ is $Q$-adapted.)
Then we have $\prD = \{ \beta_k \mid k=1, \ldots, \ell \}$, where
$$
\beta_k \seteq r_{i_1} \cdots r_{i_{k-1}}(\alpha_{i_k}) \qquad \text{ for $ k=1, \ldots, \ell$.}
$$
The reduced expression $\rxw$ gives the convex order on $\prD$ defined by $\beta_{a} > \beta_{b}$ for any $a > b$, and provides
the PBW vectors $\{ E(\beta_{k}) \}_{k=1, \ldots, \ell} $ of the negative half $U^-_{\Z[q,q^{-1}]}(\g_0)$ and its dual vectors $\{ E^*(\beta_{k}) \}_{k=1, \ldots, \ell} $ of 
the quantum unipotent coordinate ring $A_{q}(\n_0)_{\Z[q,q^{-1}]}$, where $\n_0$ is the maximal nilpotent subalgebra of $\g_0$.

Let $\{ \Cp_k \}_{k=1, \ldots, \ell} \subset R_{\g_0}\gmod$ be the \emph{cuspidal modules} associated with the reduced expression $\rxw$ (\cite{Mc15, Kato14}).
Under the isomorphism between $K(R_{\g_0}\gmod)$ and $A_{q}(\n_0)_{\Z[q,q^{-1}]}$, the cuspidal module $\Cp_k$ corresponds to the dual PBW vector $E^*(\beta_k)$ for $k=1, \ldots, \ell$.
It is known that the modules 
$$
\{  \hd ( \Cp_\ell^{\circ a_\ell} \conv \cdots \conv \Cp_1^{\circ a_1} ) \mid (a_1, \ldots, a_\ell) \in \Z_{\ge0}^\ell  \}
$$
gives a complete set of pairwise non-isomorphic simple graded $R_{\g_0}$-modules (\cite{Mc15, Kato14}).
Here $\conv$ denotes the convolution product in $R_{\g_0}\gmod$.

We now introduce the notion of cuspidal modules for quantum affine algebras.

\smallskip 

\begin{definition}
We define a sequence of simple $U_q'(\g)$-modules $\{ \cuspS_k \}_{ k\in \Z } \subset \catCO $ as follows:
\bna
\item $\cuspS_k = \F_\qQ(\Cp_k)$ for any $k=1, \ldots, \ell$, 
and we extend its definition to all $k\in\Z$ by 
\item $\cuspS_{k+\ell} = \dual( \cuspS_k )$ for any $k\in \Z$.
\ee
The modules $\cuspS_k$ $(k\in \Z)$ are called the \emph{cuspidal modules} corresponding to the Q-data $\qQ$ and the reduced expression $\rxw$.
\end{definition}

It is known that the cuspidal modules are fundamental representations 
if $\rxw$ is $Q$-adapted.
However, they may not be fundamental in general.

We have the following property.
\Prop \label{prop:unmixed}
\bnum
\item $\cuspS_k$ is a root module for any $k\in\Z$, i.e.,
$\cuspS_k$ is a real simple module satisfying \eqref{eq:root}.
\item
 $\de\bl\cuspS_j,\dual(\cuspS_k)\br=0$ for any $j<k$.
\ee
\enprop

\medskip
We define 
\begin{align} \label{Eq: ZZ}
 \ZZ \seteq \Z_{\ge0}^{\oplus \Z}
= \bigl\{ (a_k)_{k\in \Z} \in   \Z_{\ge0}^{\Z} \mid
  \sum_{k\in \Z} a_k < \infty 
\bigr\}.
\end{align}
We denote by $\prec$ the bi-lexicographic order on $\ZZ$, i.e.,   
 for any $ \bfa = (a_k)_{k\in \Z}$ and $\bfa' = (a_k')_{k\in Z}$  in $\ZZ$, $ \bfa \prec \bfa' $ if and only if the following conditions hold:
\eq\label{eq:bilexico}
&\phantom{aa}&\left\{\parbox{39ex}{\begin{enumerate}
\item there exists $r \in \Z$ such that $ a_k = a_k' $ for any $k < r$ and $ a_r < a_r'$,
\item   there exists $s \in \Z$ such that $ a_k = a_k' $ for any $k > s$ and $ a_s < a_s'$.
\end{enumerate}
}\right.
\eneq

For $ \bfa = (a_k)_{k\in \Z} \in \ZZ$, we define
\eqn
&&\sP_{Q, \rxw} (\bfa)\seteq
 \bigotimes_{k =+\infty}^{-\infty}  \cuspS_k^{\otimes a_k}
=\cdots\tens\cuspS_2^{\otimes a_2}\tens\cuspS_1^{\otimes a_1} \otimes  \cuspS_{0}^{\otimes a_{0}}
 \otimes \cuspS_{-1}^{\otimes a_{-1}}
\otimes\cuspS_{-2}^{\otimes a_{-2}}\otimes \cdots.
\eneqn
Here $\sP_{Q, \rxw} (0)$ should be understood as 
the trivial module $\one$.
We call the modules $\sP_{Q, \rxw} (\bfa)$ \emph{standard modules} with respect to the cuspidal modules $\{ \cuspS_k \}_{k\in \Z}$.  

\smallskip 

\begin{theorem} \label{Thm: main1} \label{Thm: main pbw}
\bni
\item For any $\bfa \in \ZZ$, the head of $\sP_{Q, \rxw} (\bfa)$ is simple. We denote the head by 
$$
 \sV _{Q, \rxw} (\bfa) \seteq  \hd \bl\sP_{Q, \rxw} (\bfa)\br.
$$

\item For any simple module $V \in \catCO$, there exists a unique $\bfa \in \ZZ$ such that 
$$
V \simeq  \sV _{Q, \rxw} (\bfa).
$$
\ee
Therefore the set $\{   \sV _{Q, \rxw} (\bfa) \mid \bfa \in \ZZ \}$ is a complete and irredundant set of simple modules of $\catCO$ up to isomorphisms.
\end{theorem}

\smallskip 
Indeed, (i) follows from Proposition~\ref{prop:unmixed}
(\cite{KKOP19C}).

The element $\bfa \in \ZZ$ associated with a simple module $V$ in Theorem~\ref{Thm: main1}~(ii) is called the \emph{cuspidal decomposition} of $V$ with respect to the cuspidal modules $\{ \cuspS_k \}_{k\in \Z}$,
and it is denoted by
 \begin{align} \label{Eq: cusp decomp}
\bfa_{Q, \rxw}(V) \seteq \bfa.
\end{align}

\smallskip

\begin{theorem} \label{Thm: main tri}
Let $\bfa$ be an element of $ \ZZ$. Then we have the following.
\bni
\item The simple module $ \sV _{Q, \rxw} (\bfa) $ appears only once in $\sP_{Q, \rxw} (\bfa)$.
\item If $V$ is a simple subquotient of $\sP_{Q, \rxw} (\bfa)$ which is not isomorphic to $ \sV _{Q, \rxw} (\bfa) $, then we have 
$$
\bfa_{Q, \rxw}(V) \prec \bfa.
$$
\item In the Grothendieck ring, we have 
$$
[\sP_{Q, \rxw} (\bfa)] = [\sV_{Q, \rxw} (\bfa)] + \sum_{\bfa' \prec \bfa} c(\bfa') [\sV_{Q, \rxw} (\bfa')],
$$
for some $ c(\bfa') \in \Z_{\ge 0}$.
\ee
\end{theorem}

\smallskip

\begin{remark} \label{rmk: generlization of section PBW}
In \cite{KKKO16D,KO18,OT19,OhSuh19}, the categories $\catCQ \subset \catCO$ for untwisted non simply-laced affine types and
twisted affine types are introduced, and a strong duality datum $ \ddD= \{ \Rt_i \}_{i\in J}$ for each $\catCQ$ is given.  Then, we
can obtain the duality functor
$$\F_\ddD \col R_\cmC\gmod \longrightarrow \catCQ \subset \catCO,$$
where $\cmC$ denotes the simply-laced finite Cartan matrix \emph{determined}
by the affine type of $\g$ (cf.\ \cite{KKOP20}). 
Using the functors $\F_\ddD$,
we can extend all the results in Section~\ref{sec: PBW} to an arbitrary quantum affine algebra with an arbitrary choice of $\rxw$ in type of $\cmC$.
\end{remark}

\end{document}